\documentclass[a4paper,12pt]{article}
\usepackage{amssymb,theorem}

\newtheorem{Theorem}[subsection]{Theorem}
\newtheorem{Lemma}[subsection]{Lemma}

\newtheorem{Conjecture}[subsection]{Conjecture}

\theorembodyfont{\normalfont}
\newtheorem{Remark}[subsection]{Remark}
\newtheorem{Exercise}[subsection]{Exercise}
\newtheorem{Problem}[subsection]{Problem}
\newtheorem{Example}[subsection]{Example}
\newtheorem{Definition}[subsection]{Definition}

\newcommand{\Tot}{\mathop{\mathrm{Tot}}\nolimits}
\newcommand{\End}{\mathop{\mathrm{End}}\nolimits}

\newcommand{\Mod}{\mathop{\mathsf{Mod}}\nolimits}

\newcommand{\hull}{\mathop{\mathrm{hull}}\nolimits}
\newcommand{\even}{{\mathrm{even}}}
\newcommand{\id}{{\mathrm{id}}}
\newcommand{\opp}{{\mathrm{op}}}
\newcommand{\triv}{{\mathrm{triv}}}

\newcommand{\hgt}{\mathop{\mathrm{ht}}\nolimits}
\newcommand{\gr}{\mathop{\mathrm{gr}}\nolimits}

\newcommand{\GL}{\mathop{\mathit{GL}}\nolimits}
\newcommand{\gl}{\mathop{\mathfrak{gl}}\nolimits}
\newcommand{\SL}{\mathop{\mathit{SL}}\nolimits}
\newcommand{\Ext}{\mathop{\mathrm{Ext}}\nolimits}
\newcommand{\Hom}{\mathop{\mathrm{Hom}}\nolimits}
\newcommand{\Rg}{\mathop{\mathsf{Rg}}\nolimits}

\newcommand{\Gp}{\mathop{\mathsf {Gp}}\nolimits}

\newcommand{\F}{\mathbb F}

\newcommand{\Ga}{{\mathbb G_a}}
\newcommand{\Gm}{{\mathbb G_m}}
\newcommand{\C}{{\mathsf C}}
\newcommand{\V}{{\mathcal V}}

\newcommand{\GVd}{{\mathop{\Gamma^d\mathcal V}}}
\newcommand{\GVe}{{\mathop{\Gamma^e\mathcal V}}}
\newcommand{\GVde}{{\mathop{\Gamma^{de}\mathcal V}}}
\newcommand{\Pol}{{\mathcal P}}
\newcommand{\DP}[1]{\mathop{{\mathcal D}^+{\mathcal P}_{#1}}}
\newcommand{\DbP}[1]{\mathop{{\mathcal D}^b{\mathcal P}_{#1}}}
\newcommand{\KP}{\mathop{\mathcal K^+\mathcal P}}
\newcommand{\KI}{\mathop{\mathcal K^+\mathcal I}}
\newcommand{\KbI}{\mathop{\mathcal K^b\mathcal I}}
\newcommand{\KrA}{{{\mathbf K}^{\mathbf r}_A}}
\newcommand{\Kr}{{{\mathbf K}^{\mathbf r}}}
\newcommand{\Char}{\mathop{\mathrm {Char}}}

\newcommand{\CC}{{\mathbb C}}
\newcommand{\cS}{{\mathfrak S}}

\newcommand{\cL}{{\cal L}}
\newcommand{\qed}{\unskip\nobreak\hfill\hbox{ $\Box$}}

\title{Lectures on bifunctors and finite generation of rational cohomology algebras}
\author{Wilberd van der Kallen}
\date{}
\begin{document}\sloppy

\maketitle

\begin{abstract}
This text is an updated version of material used for a course at Universit\'e de Nantes, part of 
`Functor homology and applications', April 23-27, 2012.
The  proof \cite{Touze 2010}, \cite{TvdK} by Touz\'e of my conjecture on cohomological finite
generation (CFG) has been one of the successes of functor homology. We will not treat this proof in any detail. 
Instead we will focus on a formality conjecture of Cha\l upnik
and discuss ingredients of a second generation proof \cite{Touze 2013} of the existence of the universal classes of 
Touz\'e.
\end{abstract}
\tableofcontents
\section{The CFG theorem}
In its most basic form the CFG theorem of \cite{TvdK} reads
\begin{Theorem}[Cohomological Finite Generation]
 Let $G$ be a reductive algebraic group over an algebraically closed field $k$, and let $A$ be a finitely generated 
commutative $k$-algebra  on which $G$ acts algebraically via algebra automorphisms.
Then the cohomology algebra $H^*(G,A)$ is a finitely generated graded
$k$-algebra.
\end{Theorem}

An essential ingredient in the proof of this theorem is the existence of certain universal 
cohomology classes. They were constructed by Touz\'e in \cite{Touze 2010}. We  will discuss methods used in the new
construction \cite{Touze 2013} of
these classes. 

\section{Some history}

Let us give some background. First there is \emph{invariant theory} \cite{Borel history}, \cite{Grosshans book}, 
\cite{Springer inv book}.
Classical invariant theory looked at the following situation. (We will give a very biased description, full of anachronisms.)
Say the algebraic Lie group $G(\CC):=\SL_n(\CC)$ acts on a finite dimensional complex vector space
$V$ with dual $V^\vee$. 
Then it also acts on the symmetric algebra $A=S_\CC^*(V^\vee)$ of polynomial maps from $V$ to $\CC$.
One is interested in the subalgebra $A^{G(\CC)}$ of elements fixed by $G(\CC)$. It is called the subalgebra of \emph{invariants}.
More generally, if $W$ is another finite dimensional complex vector space on which $G(\CC)$ acts, then $ W\otimes_\CC A$
encodes the polynomial maps from $V$ to $W$. The subspace $ (W\otimes_\CC A)^G$ of
fixed points or invariants in $ W\otimes_\CC A$ corresponds with the
equivariant polynomial maps from $V$ to $W$.
This subspace of invariants is a module over the algebra of invariants $A^{G(\CC)}$.
When $n=2$ and $V$ is irreducible Gordan (1868) showed  in a constructive manner that 
$A^{G(\CC)}$ is a finitely generated algebra 
\cite{Gordan}.
Our $V$ corresponds with his  `binary forms of degree $d$', with $d=\dim V-1$.
Hilbert (1890) generalized Gordan's theorem  nonconstructively to arbitrary $n$ 
and---encouraged by an incorrect claim of Maurer---%
asked in his 14-th problem to prove that this finite generation of invariants is a very general fact about actions of
algebraic Lie groups on domains of finite type
over $\CC$. A counterexample of Nagata (1959) showed this was too optimistic, but by then it was understood that finite 
generation of invariants holds for compact connected real Lie groups (cf. Hurwitz 1897)
as well as for their complexifications, also known as the connected
reductive complex algebraic Lie groups (Weyl 1926). 
Finite groups have been treated by Emmy Noether (1926) \cite{Noether}, so connectedness may be dropped.
(Algebraic Lie groups have finitely many connected components.)

Mumford (1965) needed finite generation of invariants for reductive
algebraic groups over fields of arbitrary
characteristic in order to construct moduli spaces. 
In his book Geometric Invariant Theory \cite{GIT} he introduced a condition, often referred to as \emph{geometric reductivity}, 
that he conjectured to be true
for reductive algebraic groups and that he conjectured to imply finite generation of invariants. These conjectures
were confirmed by Haboush (1975) \cite{Haboush}
and Nagata (1964) \cite{Nagata} respectively. Nagata treated any algebra of finite type over the base field, not just domains.
We adopt this generality. It rather changes the problem of finite generation of invariants.
For instance, counterexamples to finite generation of invariants are now easy to find already when the Lie group $G(\CC)$
is  $\CC$ with addition as operation. (See Exercise \ref{not reductive}.)

[We now understand that over an arbitrary commutative noetherian
base ring the right counterpart
of Mumford's geometric reductivity is not the geometric reductivity of Seshadri (1977) but the \emph{power reductivity} of
Franjou and van der Kallen (2010) \cite{FvdK}, which is actually equivalent to the finite generation property.]

Let us say that $G$ satisfies property (FG) if, whenever $G$ acts on a commutative algebra of $A$ finite type over $k$,
the ring of invariants $A^G$ is also finitely generated over $k$. So then the theorem of Haboush and Nagata says that 
connected reductive algebraic groups over a field have property (FG). Of course the action of $G$ on $A$ should be consistent 
with the nature of $G$ and $A$ respectively. Thus if $G$ is an algebraic group, then the action should be algebraic and
 the multiplication 
map $A\otimes_kA\to A$ should be equivariant.

We will be interested in the cohomology algebra $H^*(G,A)$ of a geometrically reductive group $G$ acting on a commutative
algebra $A$ of finite type over a base field $k$. Or, more generally, a power reductive affine flat algebraic group
scheme $G$ acting on  a commutative algebra $A$ of
finite type  over a noetherian commutative base ring $k$.
Observe that $H^0(G,A)$ is just the algebra of invariants $A^G$, which we know to be finitely generated.
The $H^i(G,{-})$ are the right derived functors of the functor $V\mapsto V^G$.

My conjecture was that the full algebra $H^*(G,A)$ is finitely generated when $k$ is field and $G$ is a
geometrically reductive group (or group scheme).
Let us say that $G$ satisfies the cohomological finite generation property (CFG) if,
 whenever $G$ acts on a commutative algebra $A$ of finite type over $k$,
the cohomology algebra $H^*(G,A)$ is also finitely generated over $k$.
So my conjecture was that if the base ring $k$ is a field and an affine algebraic group (or group scheme) $G$ over $k$
 satisfies property (FG) then it actually satisfies the stronger property (CFG).
This was proved by Touz\'e \cite{Touze 2010}, 
by constructing classes $c[m]$ in Ext groups in the category of \emph{strict polynomial bifunctors} of
Franjou and Friedlander \cite{Franjou-Friedlander}. If the base field has characteristic zero then there is little to do,
because then (FG) implies that $H^{>0}(G,A)$ vanishes.

One may ask  if (CFG) also holds
 when the base ring is not a field but just noetherian and $G=\GL_n$ say.
This question is still open for $n\geq3$. But see \cite{vdK family}.

We are not aware of striking applications of
 the general (CFG) theorem, but investigating the (CFG) conjecture has led to new insights \cite{vdK gfd}.
The conjecture also fits into a long story where special cases have been
very useful. The case of a finite group was treated by Evens (1961) \cite{Evens}
and this has been the starting point for the theory
of \emph{support varieties} \cite[Chapter 5]{Benson II}. In this theory one exploits a connection between the 
rate of growth of a minimal 
projective resolution and 
the dimension of a `support variety', which is a subvariety of the spectrum of $H^\even(G,k)$.
The case of finite group schemes over a field (these are group schemes whose coordinate ring is a finite dimensional 
vector space) turned out to be `surprisingly elusive'. It was finally settled by Friedlander and Suslin (1997) 
\cite{Friedlander-Suslin}.
For this they had to invent \emph{strict polynomial functors} and compute with certain Ext groups in the category of strict 
polynomial functors.
Again their result was crucial for developing a theory of support varieties, now for finite group schemes.

As $H^{>0}(G,k)$ vanishes for reductive $G$, there is no obvious theory of support varieties for reductive $G$. 

\begin{Exercise}[Additive group is not reductive]\label{not reductive}
 Let $G=\CC$ with addition as group operation. Make $G$ act on $M=\CC^2$ by $x\cdot(a,b)=(a+xb,b)$. Projection
onto the second factor of $\CC^2$ defines a  surjective
equivariant linear map $M\to\CC$ with $G$ acting trivially
on the target. It induces a map of symmetric algebras $S_\CC^*(M)\to S_\CC^*(\CC)$. View
$S_\CC^*(\CC)$ as an $S_\CC^*(M)$-module. Show that the algebra of invariants in the finite type $\CC$-algebra
$S_{S_\CC^*(M)}^*(S_\CC^*(\CC))$ is not finitely generated. Hint: Exploit the trigrading.
\end{Exercise}

\emph{Reductivity} can be thought of as what one needs to avoid this example and its relatives. 
Reductivity of an affine algebraic group $G$ over an algebraically closed field
forbids that the connected component of the identity of $G$ 
(for the Zariski topology) has a
normal algebraic subgroup isomorphic to the additive group underlying a nonzero vector space.
Originally reductivity referred to representations being completely reducible, but this meaning was abandoned in order 
to include groups over fields of positive characteristic that look pretty much like reductive groups over $\CC$.
For example $\GL_n$ is reductive, but when the ground field has positive characteristic, $\GL_n$ has representations 
that are not completely reducible. Indeed in positive characteristic the category of representations of $\GL_n$ has  
interesting $\Ext$ groups and this is our subject.

\section{Some basic notions, notations and facts for group schemes}
Let us now assume less familiarity with algebraic groups or group schemes.

\subsection{Rings and algebras}
Every ring has a unit and ring homomorphisms are unitary.
Our \emph{base ring} $k$ is commutative noetherian and most of the time a field of characteristic $p>0$, in fact just $\F_p$.
Let $\Rg_k$ denote the category of commutative $k$-algebras. An object $R$ of $\Rg_k$ is a commutative ring together with
a homomorphism $k\to R$. We write $R\in \Rg_k$ to indicate that $R$ is an object of  $ \Rg_k$. 
The same convention will be used for other categories. When $\C$ is a category, $\C^\opp$ denotes the opposite category.
Let $\Gp$  be the category of groups. 

\subsection{Group schemes}
A functor $G:\Rg_k\to \Gp$ is called an affine flat algebraic \emph{group
scheme} over $k$ if $G$ is \emph{representable} \cite[1.2]{Waterhouse book}, \cite{MacLane working} 
by a flat $k$-algebra of finite type, which is then known as the \emph{coordinate ring}
$k[G]$ of $G$  \cite{Demazure-Gabriel}, \cite{Jantzen book}, \cite{Waterhouse book}. 
Recall that this means that for every $R$ in $\Rg_k$ one is given a bijection between $\Hom_{\Rg_k}(k[G],R)$ and 
$G(R)$, thus providing $\Hom_{\Rg_k}(k[G],R)$ with a group structure, functorial in $R$. 
In particular one has the unit element $\epsilon:k[G]\to
k$ in the group $G(k)\cong \Hom_{\Rg_k}(k[G],k)$. This $\epsilon$ is also known as the \emph{augmentation map} of $k[G]$.
In the group $\Hom_{\Rg_k}(k[G],k[G]\otimes_k k[G])$ one has the elements $x:f\mapsto f\otimes 1$ and 
$y:f\mapsto 1\otimes f$ with 
product $xy$ known as the \emph{comultiplication} $\Delta_G:k[G]\to k[G]\otimes_k k[G]$. These maps $\epsilon$,
$\Delta_G$ make $k[G]$ into a \emph{Hopf algebra} \cite[1.4]{Waterhouse book}. 
(There is also an \emph{antipode}.)
If $g$, $h\in \Hom_{\Rg_k}(k[G],R)$
then $gh$ in $G(R)$ is just $m_R\circ( g\otimes h)\circ \Delta_G$, where $m_R:R\otimes_kR\to R$ is the multiplication map of $R$.

\subsection{$G$-modules}
We will be working in the category $\Mod_G$ of $G$-modules. A $G$-module or \emph{representation} of $G$ is simply a 
\emph{comodule} \cite[3.2]{Waterhouse book}
for the Hopf algebra $k[G]$. In functorial language this means that one is given a $k$-module $V$
with an action of $G(R)$ on $V\otimes_k R$ by $R$-linear endomorphisms, functorially in the commutative $k$-algebra $R$.
In particular, the identity map $k[G]\to k[G]$ viewed as an element of $G(k[G])$ acts by a $k[G]$-linear map 
$V\otimes_kk[G]\to  V\otimes_kk[G]$ and the composite of this $k[G]$-linear map with $v\mapsto v\otimes 1$ is the 
\emph{comultiplication} 
$\Delta_V:V\to V\otimes_kk[G]$ defining the comodule structure of $V$.
If $g\in \Hom_{\Rg_k}(k[G],R)\cong G(R)$, then it acts on $V\otimes_kR$ as $v\otimes r\mapsto (\id\otimes g)(\Delta_V(v))r$.
The category $\Mod_G$ has useful properties only under the assumption that $G$ is flat over $k$. 
That is why we always make this assumption. 
Flatness is of course automatic when $k$ is a field. Geometers should be warned that it is a mistake to restrict attention
to representations that are representable. So while our group functors are schemes, our representations need not be.
For instance, in the (CFG) conjecture finite dimensional algebras $A$ are of less interest. 
And if $A$ is infinite dimensional as a vector space
then as a representation it is no scheme.

\subsection{Invariants}One may define the submodule $V^G$ of fixed
vectors or \emph{invariants} of a representation $V$ and get a natural isomorphism 
$\Hom_{\Mod_G}(k,V)\cong V^G$, where $k$ also stands 
for the representation $k^\triv$ with underlying module $k$ and trivial $G$ action. 

\subsection{Cohomology of $G$-modules}The category $\Mod_G$ is abelian with enough injectives. 
We write $\Hom_G$ for $\Hom_{\Mod_G}$ and $\Ext_G$ for $\Ext_{\Mod_G}$. 
Cohomology is simply defined as follows: $$H^i(G,V):=\Ext^i_G(k,V).$$
It may be computed \cite[I 4.14--4.16]{Jantzen book} as the cohomology of the Hochschild complex 
$C^\bullet(V)=(V\otimes_k C^\bullet(k[G]))^G$. 
There is a \emph{differential graded algebra} (=DGA) structure on $C^\bullet(k[G])=k[G]^{\otimes(\bullet+1)}$. 
Let $R\in \Rg_k$ be provided with an action of $G$. So $R$ is a $G$-module and the multiplication $R\otimes_kR\to R$ is a
$G$-module map. 
If $u\in C^r(G,R)$ and $v\in C^s(G,R)$, then $u\cup v$
is defined in simplified notation by 
$$(u\cup v)(g_1,\ldots g_{r+s})=u(g_1,\ldots,g_r).{}^{g_1\cdots g_r} v(g_{r+1},\ldots ,g_{r+s}),$$
where ${}^gr$ denotes the image of $r\in R$ under the action of $g$.
With this cup product $C^*(G,R)$ is a differential graded algebra.

\begin{Remark}
 We have followed \cite{Jantzen book} in that we have used inhomogenous cochains, although for $C^\bullet(k[G])$ homogeneous cochains
might be more natural. Thus one could take as alternative starting point a differential graded algebra
 $C_{\mathrm {hom}}^\bullet(k[G])$  with 
$C_{\mathrm {hom}}^i(k[G])=k[G]^{\otimes(i+1)}$
and differential $d$ as suggested by $(df)(g_0,g_1,g_2)=f(g_1,g_2)-f(g_0,g_2)+f(g_0,g_1)$. View $C_{\mathrm {hom}}^i(k[G])$ as 
$G$-module through left translation as in ${}^gf(g_0,\cdots,g_i)=f(g^{-1}g_0,\cdots,g^{-1}g_i)$. 
Then $H^i(G,V)$ may be computed as the cohomology of 
$(V\otimes_k C_{\mathrm {hom}}^\bullet(k[G]))^G$.
\end{Remark}

\subsection{Symmetric and divided powers}
For simplicity let $k$ be a field. If $V$ is a finite dimensional vector space and $n\geq1$, we have an action of
the \emph{symmetric group} $\cS_n$ on $V^{\otimes n}$ and the $n$-th \emph{symmetric power} $S^n(V)$ is the module of
\emph{coinvariants} \cite[II 2]{Brown} $(V^{\otimes n})_{\cS_n}=H_0(\cS_n,V^{\otimes n})$ for this action.
Dually the $n$-th \emph{divided power} $\Gamma^n(V)$ is the module of invariants $(V^{\otimes n})^{\cS_n}$ \cite{Eisenbud book},
\cite{Roby}.
One has $\Gamma^n(V)^\vee\cong S^n(V^\vee)$.

Both $S^*$ and $\Gamma^*$ are \emph{exponential functors}.
That is, one has $$S^n(V\oplus W)=\bigoplus_{i=0}^n S^i(V)\otimes_kS^{n-i}(W)$$
and similarly $$\Gamma^n(V\oplus W)=\bigoplus_{i=0}^n \Gamma^i(V)\otimes_k\Gamma^{n-i}(W).$$

\subsection{Tori} 
 A very important example of an algebraic group scheme is the \emph{multiplicative group} $\Gm$. It associates to $R$ its group of invertible elements $R^*$.
The coordinate ring $k[\Gm]$ is the Laurent polynomial ring $k[X,X^{-1}]$. 
Any $\Gm$-module $V$ is a direct sum of \emph{weight spaces}
$V_i$ on which $\Delta_V$ equals $v\mapsto v\otimes X^i$. Weight spaces are nonzero by definition.

\begin{Exercise}
 Prove this decomposition into weight spaces. 
Rewrite $k[X,X^{-1}]\otimes_kk[X,X^{-1}]$ as $k[X,X^{-1},Y,Y^{-1}]$ where  $X\otimes 1$ is written as $X$ and
$1\otimes X$ as $Y$, so that $\Delta_{\Gm}X=XY$.
Use that if $\Delta_Vv=\sum_i\pi_i(v)X^i$,  then
$\sum_{i}\pi_i(v)(XY)^i=\sum_{i,j}\pi_j(\pi_i(v))X^iY^j$.
\end{Exercise}

More generally the direct product $T$ of $r$ copies of $\Gm$, known as a \emph{torus $T$ of rank $r$} has
as coordinate ring the Laurent polynomial ring in $r$ variables $k[X_1,X_1^{-1},\cdots,X_r,X_r^{-1}]$.
Again any $T$-module $V$ is a direct sum of nonzero \emph{weight spaces} $V_\lambda$ where now the 
weight $\lambda$ is an $r$-tuple of
integers and $\Delta_V$ restricts to $v\mapsto v\otimes X_1^{\lambda_1}\cdots X_r^{\lambda_r}$ on $V_\lambda$.
So a weight space is spanned by simultaneous eigenvectors with common eigenvalues and every $T$-module is diagonalizable.
The invariants in a $T$-module are the elements of weight zero. Taking invariants is exact on $\Mod_T$ and
$H^{>0}(T,V)$ always vanishes.

\subsection{The additive group}
The group scheme $\Ga$ sends a $k$-algebra $R$ to the underlying additive group. The coordinate ring of $\Ga$ is $k[X]$ with
$\Delta_{\Ga}(X)=X\otimes1+1\otimes X$.
Recall that the additive group is not reductive. It has no property (FG). (Redo exercise \ref{not reductive} with $k$ 
replacing $\CC$.)
If $k$ is a field of characteristic $p>0$ then $H^1(\Ga,k)$ is already infinite
dimensional, so even with such small coefficient module the cohomology explodes.
Thus cohomological finite generation is definitely tied with reductivity.

\subsection{General linear group}
Let $n\geq1$. The group scheme $\GL_n$ associates to $R$ the group $\GL_n(R)$ of $n$ by $n$ matrices with entries in $R$
and with invertible determinant. Its coordinate ring $k[\GL_n]$ is $k[M_n][1/{\det}]$, where $k[M_n]$, 
also known as the coordinate
ring of the \emph{monoid} of $n$ by $n$ matrices, is the polynomial ring $k[X_{11},X_{12},\cdots,X_{nn}]$
in $n^2$ variables $X_{11},X_{12},\cdots,X_{nn}$ and
$\det$ is the determinant of the matrix $(X_{ij})$. A ring homomorphism $\phi:k[M_n]\to R$ corresponds with the matrix
$(\phi(X_{ij}))$ and $\phi$ extends to $k[\GL_n]$ if and only if this matrix is invertible. One sees that 
indeed $\Hom_{\Rg_k}(k[\GL_n],R)\cong \GL_n(R)$. If $n=1$ we are back at $\Gm$, but as soon as $n\geq2$ the representation
theory becomes much more interesting. In fact there is a lemma (cf.\ \cite[Lemma 1.7]{TvdK}) telling that for proving  
my (CFG) conjecture over a field $k$ 
it suffices to show that the reductive group scheme $G=\GL_n$ has (CFG), in particular for large $n$.
The lemma explains why the homological algebra of strict
polynomial bifunctors becomes so relevant: As we will see, it encodes  what happens
 to 
$H^\bullet({\GL_n},{V_n})$ as $n$ becomes large, for a certain kind of coefficients
$V_n$.

\subsection{Polynomial representations}
Let $k$ be a field until further notice. 
One calls a finite dimensional representation $V$ of $\GL_n$ a \emph{polynomial representation} 
if the action
is given by polynomials, meaning that $\Delta_V$ factors trough the embedding $V\otimes_kk[M_n]\to V\otimes_kk[\GL_n]$.
And one calls it homogeneous of \emph{degree} $d$ if moreover $\Delta_V$ lands in $V\otimes_kk[M_n]_d$, where $k[M_n]_d$
consists of polynomials homogeneous of total degree $d$. If one lets $\Gm$ act on $k[M_n]$ by algebra automorphisms
giving the variables $X_{ij}$ weight one and $k$ weight zero, then $k[M_n]_d$ is just the weight space of weight $d$.
Polynomial representations were studied by Schur in his thesis (1901).
The \emph{Schur algebra} $S_k(n,d)$ can be described as 
$\Gamma^d(\End_k(k^n))$ 
with multiplication obtained by restricting the usual algebra structure on $\End_k(k^n)^{\otimes d}$ given
by $(f_1\otimes\cdots\otimes f_d)(g_1\otimes\cdots\otimes g_d)=f_1g_1\otimes\cdots\otimes f_dg_d$.
The category of finitely generated left $S_k(n,d)$-modules is equivalent to the category of finite dimensional
polynomial representations of degree $d$ of $\GL_n$ \cite[\S 3]{Friedlander-Suslin}.

\subsection{Frobenius twist of a representation}
Let  $p$ be a prime number and $k=\F_p$.
The group scheme $\GL_n$ admits a Frobenius homomorphism $F:\GL_n\to\GL_n$ that sends a matrix $(a_{ij})\in \GL_n(R)$
to $(a_{ij}^p)$. If $V$ is a representation of $\GL_n$ then one gets a new representation $V^{(1)}$, called the 
\emph{Frobenius twist},
by precomposing with $F$.
If $V$ is a polynomial representation of degree $d$ then $V^{(1)}$ has degree $pd$.
One may also twist $r$ times and obtain $V^{(r)}$. We do not reserve the notation $F$ for Frobenius, but $V^{(r)}$
will always indicate an $r$-fold Frobenius twist.

\begin{Exercise}
We keep $k=\F_p$. Let $V$ be a finite dimensional representation of $\GL_n$. Choose a basis in $V$. 
The action of $g\in\GL_n(R)$ on
$V\otimes_kR$ is given with respect to the chosen basis by a matrix $(g_{ij})$ with entries in $R$. 
Show that the action on $V^{(r)}\otimes_kR$ is given by
the matrix $(g_{ij}^{p^r})$. In other words, when the base field is $\F_p$ one may confuse precomposition by Frobenius with 
postcomposition. For larger ground fields one would have to be more careful.
\end{Exercise}

\section{Some basic notions, notations and facts for functors}

\subsection{Strict polynomial functors}
Let $\V_k$ be the $k$-linear category of finite dimensional vector spaces over a field $k$.
The category $\GVd_k$, often written $\GVd$, generalizes the {Schur algebras} as follows.
Its objects are finite dimensional vector spaces over $k$, but $\Hom_{\GVd}(V,W)=\Gamma^d(\Hom_k(V,W))$.
The composition is similar to the one in a Schur algebra. We could call $\GVd_k$ the \emph{Schur category}.
The category of \emph{strict polynomial functors of degree $d$} is now defined, following 
the exposition of Pirashvili 
\cite{Pirashvili}, \cite{Pirashvili B}, 
as the category of $k$-linear functors $\GVd\to\V_k$. The reason for the word \emph{strict} is simply that
the terminology \emph{polynomial functor} already means something. There is an obvious functor $\iota^d$
from $\V_k$ to $\GVd$.
It sends $V\in\V_k$ to $V\in\GVd$ and $f\in \Hom_k(V,W)$ to $f^{\otimes d}$. This is not $k$-linear when $d>1$.
If $F\in \Pol_d$, let us try to understand the composite map $\Hom_k(V,W)\to \Hom_k(FV,FW)$.
The map $\Hom_\GVd(V,W)\to \Hom_k(FV,FW)$ is $k$-linear and is thus given by an element $\psi$ of the space
$\Hom_k(\Gamma^d(\Hom_k(V,W)),\Hom_k(FV,FW))\cong \Hom_k(FV,FW)\otimes S^d(\Hom_k(V,W)^\vee)$ which also encodes the
polynomial maps from $\Hom_k(V,W)$ to $  \Hom_k(FV,FW)$ that are homogeneous of degree $d$.
One checks that the composite map $\Hom_k(V,W)\to \Hom_k(FV,FW)$ is the polynomial map of degree $d$ encoded by $\psi$. 
This explains why $F$ is called
a (strict) polynomial functor of degree $d$.

\begin{Remark}
 The original definition of Friedlander and Suslin did not use $\GVd$, but just defined strict polynomial functors 
of degree $d$ as
functors $F:\V_k\to \V_k$ enriched with elements $\phi_{V,W}$ in $\Hom_k(FV,FW)\otimes S^d(\Hom_k(V,W)^\vee)$
that satisfy appropriate conditions, like the condition that the polynomial map $\Hom_k(V,W)\to \Hom_k(FV,FW)$
encoded by $\phi_{V,W}$ agrees with $F$. 
That is more intuitive, but the definition by means of
$\GVd$ is concise and  has its own advantages. In fact one may view $F:\GVd\to \V_k$ as exactly the enrichment that 
Friedlander and Suslin need to add to the composite functor $F\iota^d$.
One should use both points of view. They are equivalent \cite{Pirashvili}.
We will secretly think in terms of the Friedlander and Suslin setting when that is more convenient.
\end{Remark}

\subsection{Some examples of strict polynomial functors}
The functor $F=\otimes^d$ maps $V\in\GVd$ to $V^{\otimes d}$. If $f\in \Hom_\GVd(V,W)$, view $f$ as an element
of $\Hom_k(V,W)^{\otimes d}$ and define $Ff:FV\to FW$ by means of the pairing  $\Hom_k(V,W)^{\otimes d}\times
V^{\otimes d}\to W^{\otimes d}$ which maps the pair $(f_1\otimes\cdots\otimes f_d,v_1\otimes\cdots\otimes v_d)$
to $f_1(v_1)\otimes\cdots\otimes f_d(v_d)$.

The functor $\Gamma^d$ is the subfunctor of $\otimes^d$ with value $\Gamma^d(V)$
on $V\in\GVd$.

The functor $S^d$ is the quotient functor of $\otimes^d$ with value $S^d(V)$
on $V\in\GVd$.

If $F\in \Pol_d$, then its \emph{Kuhn dual} $F^\#$ is defined as $DFD$, where $DV=V^\vee$ is the contravariant functor
on $\V_k$ or ${\Gamma^d}{\V_k}$ sending $V$ to its $k$-linear dual $V^\vee$. Thus $S^{d\#}=\Gamma^d$.

If $k$ has characteristic $p>0$, then the $r$-th \emph{Frobenius twist functor} $I^{(r)}\in\Pol_{p^r}$ is the subfunctor 
of $S^{p^r}$ such that the vector space
$I^{(r)}V$ is
generated by the $v^{p^r}\in S^{p^r}V$. Note that every element of $I^{(r)}V$ is actually of the form $v^{p^r}$ if $k=\F_p$.

\subsection{Polynomial representations from functors }
If $F\in\Pol_d$ then $F(k^n)$ is a polynomial representation of degree $d$ of $\GL_n$.
The comodule structure is obtained from the homomorphism $\Hom_\GVd(k^n,k^n)\to \Hom(F(k^n),F(k^n))$ by means of the isomorphism 
$\Hom_k(\Hom_\GVd(k^n,k^n),\Hom(F(k^n),F(k^n)))\cong \Hom(F(k^n),F(k^n))\otimes_kS^d(\Hom_k(k^n,k^n)^\vee)$.

Friedlander and Suslin showed \cite[\S 3]{Friedlander-Suslin} that if $n\geq d$ this actually 
provides an equivalence of categories, preserving $\Ext$ groups \cite[Cor 3.12.1]{Friedlander-Suslin}, between $\Pol_d$ and 
the category of finite dimensional polynomial representations of degree $d$
of $\GL_n$. So again there is another way to look at $\Pol_d$
and we secretly think in terms of polynomial representations when we find that more convenient. 

\begin{Exercise}[Polarization]
 If $k$ is a finite field and $V\in\V_k$, then $V$ is a finite set.
If $\dim V>1$ and $d$ is large, then the dimension of $\Gamma^dV$ exceeds 
the number of elements of $V$ so that $\Gamma^dV$ is certainly
not spanned by elements of the form $v^{\otimes d}$. On the other hand Friedlander and Suslin show that $\Gamma^dV$ is
spanned by such elements if $k$ is big enough, when keeping $d$ and $\dim V$ fixed. So as long as one uses constructions
that are compatible with base change one may think of $\Gamma^dV$ as spanned by the $v^{\otimes d}$.

Let $V=k^n$. Show that the $v^{\otimes d}$
generate $V$ as a $\GL_n$-module. Hint: Let $T$ be the group scheme of diagonal matrices in $\GL_n$. 
Show that the weight spaces 
of $T$ in $\Gamma^dV$ are one dimensional. Any $\GL_n$-submodule must be a $T$-submodule, hence a sum of weight spaces.
Now compute the weight decomposition of $v^{\otimes d}$ for $v\in V$.
\end{Exercise}

\subsection{Composition of strict polynomial functors}
If $F\in\Pol_d$, $G\in \Pol_e$ we wish to define their composite $F\circ G\in \Pol_{de}$. 
Associated to $F$ one has the functor $F\iota^d:\V_k\to\V_k$ and associated to $G$ one has $G\iota^e:\V_k\to\V_k$.
We want $F\circ G$ to correspond with the composite of $F\iota^d$ and $G\iota^e$.
For $V\in \GVde$ one puts $(F\circ G)V=F(GV)$. For $f\in \Hom_k(V,W)$ we want  that $(F\circ G)f^{\otimes de}$ 
equals
$F(G(f^{\otimes e})^{\otimes d})$. 
Thus let $\phi:\Hom_{\GVe}(V,W)\to \Hom_k(GV,GW)$ be given by $G$
and observe that the restriction of $\Gamma^d\phi:\Gamma^d\Hom_{\GVe}(V,W)\to \Gamma^d\Hom_k(GV,GW)$ to $\Gamma^{de}\Hom_k(V,W)$
lands in the source of the map $\Hom_{\GVd}(GV,GW)\to \Hom_k(FGV,FGW)$. 
\begin{Exercise}
 Finish the definition and check all details.
\end{Exercise}

In particular, the composite $F\circ I^{(r)}$ is called the $r$-th \emph{Frobenius twist $F^{(r)}$ of the functor} $F\in\Pol_d$.
Recall that if $n\geq d$ we have an equivalence of categories, between $\Pol_d$ and 
the category of finite dimensional polynomial representations of $\GL_n$ of degree $d$. Take $k=\F_p$
for simplicity. Now check that the notion of Frobenius twist on the strict polynomial side agrees with the notion of Frobenius twist
for representations.

\subsection{Untwist}\label{untwist coh}
For $F$, $G\in \Pol_d$ and $r\ge1$ we have $\Hom_{\Pol_d}(F,G)\cong\Hom_{\Pol_{dp^r}}(F^{(r)},G^{(r)})$ by 
\cite[Lemma 2.2]{Touze 2012}.
So to construct a morphism in $\Pol_d$ one may twist first. 
This was well known in the context of representations of $\GL_n$, but the proof we know  there  involves 
fppf sheaves
 \cite[I 9.5; I 6.3]{Jantzen book}.

On $\Ext_{\GL_n}$ groups Frobenius twist gives injective maps \cite[II 10.14]{Jantzen book}, but often no isomorphisms.
Compare the formality conjecture below.
In view of the connection between $\Ext_{\GL_n}$ groups and $\Ext_{\Pol_d}$ groups we may also state this twist injectivity
as 

\begin{Theorem}[Twist Injectivity]\label{twist injectivity}
Let $F,G\in \Pol_d$.
Precomposition by $I^{(1)}$ induces an injective map $\Ext^i_{\Pol_d}(F,G)\to\Ext^i_{\Pol_{dp}}(F^{(1)},G^{(1)})$
for every $i\geq0$. 
\end{Theorem}

\subsection{Parametrized functors}
If  $V\in\GVd$ define the functor  $V\otimes_k^{\Gamma^d}{-}:\GVd\to\GVd$ by sending an object $W$ to $V\otimes_kW$
and a morphism $f\in \Hom_\GVd(W,Z)$ to its image under 
$\Gamma^d(\phi):\Gamma^d\Hom_k(W,Z)\to \Gamma^d\Hom_k(V\otimes_kW,V\otimes_kZ)$ where $\phi:g\mapsto \id_V\otimes g$.
One checks that $V\otimes_k^{\Gamma^d}{-}$ is functorial in $V$.

If $F\in \Pol_d$, $V\in\GVd$, then $F_V$ denotes the composite $F(V\otimes_k^{\Gamma^d}{-})$, $F_VW=F(V\otimes_k^{\Gamma^d}W)$.
It is covariantly functorial in $V$, which is why we use a subscript.
Dually, $F^V$ denotes $F(\Hom_k(V,{-}))=((F^\#)_V)^\#$. It is contravariantly functorial in $V$, 
which is why we use a superscript. Notice that we did not decorate $\Hom_k$ with $\Gamma^d$ like we did with $\otimes_k$. 
We leave that
to the reader.

For example, $\Gamma^{dV}W=\Gamma^d(\Hom_k(V,W))=\Hom_\GVd(V,W)$, so that the Yoneda lemma \cite[1.3]{Waterhouse book},
\cite{MacLane working} gives
$$\Hom_{\Pol_d}(\Gamma^{dV},F)\cong FV.$$ As $FV$ is exact in $F$, it follows that $\Gamma^{dV}$ is projective in $\Pol_d$.
Dually $S^d_V=\Gamma^{dV\#}$ is injective in $\Pol_d$ and $$\Hom_{\Pol_d}(F,S^d_V)\cong F^\#(V).$$

\subsection{An adjunction}For $F$, $G\in\Pol_d$ we have 
$$\Hom_{\Pol_d}(F^V,G)\cong\Hom_{\Pol_d}(F,G_V)$$ in $\Pol_d$. So $F\mapsto F^V$ has 
\emph{right adjoint} \cite{MacLane working}
$G\mapsto G_V$. 
Indeed the standard map
$V\otimes_k\Hom_k(V,W)\to W$ in $\V_k$ 
induces a morphism $V\otimes_k^{\Gamma^d}\Hom_k(V,W)\to W$ in 
$\GVd$ so that if $\phi:F\to G_V$ one gets a map $F^VW=F\Hom_k(V,W)\to G_V\Hom_k(V,W)\to
GW$, functorial in $G$.
If $G=S^{d}_Z$ then $\Hom_{\Pol_d}(F,G_V)\to\Hom_{\Pol_d}(F^V,G)$ becomes the isomorphism 
$F^\#(Z\otimes_k^{\Gamma^d}V)\to F^\#(V\otimes_k^{\Gamma^d}Z)$.
As $\Hom_{\Pol_d}({-},{-})$ is left exact, the result follows from this and functoriality in $G$.

\subsection{Coresolutions}
If $\dim V\geq d$ then $\Gamma^{dV}$ forms a \emph{projective generator} \cite{MacLane working} of 
$\Pol_d$ and $S^d_V$ an injective cogenerator.
Say $V=k^n$ with $n\geq d$ and let $G=\GL_n$ again. One may also write $G=\GL_V$. 
For $F\in\Pol_d$ we have $FW\cong \Hom_{\Pol_d}(F^\#,S^d_W)\cong\Hom_G(F^\#V,S^d_WV)\hookrightarrow\Hom_k(F^\#V,S^d_WV)\cong
\Hom_k(F^\#V,S^d_VW)$, functorially in $W$, so that $$F\hookrightarrow\Hom_k(F^\#V,S^d_V).$$ And $\Hom_k(F^\#V,S^d_V)$
is just a direct sum of $\dim F^\#V$ copies of $S^d_V$, so it is injective and we conclude that $\Pol_d$ has enough injectives.
Therefore we know now how to build injective coresolutions consisting of direct sums of copies of
$S^d_V$. As $\End_{\Pol_d}(S^d_V)\cong
\End_\GVd(V)$ we also have a grip on the differentials in these coresolutions.

So far we discussed coresolving an object of $\Pol_d$.
We also want to coresolve cochain complexes.
When we speak of a \emph{cochain complex} $C^\bullet$ we do not assume $C^i$ to vanish for $i<0$.
When $f$ is a cochain map, we may use the symbol $\hookrightarrow$ to indicate that is an injective cochain map.
If $$C^\bullet  = 
\cdots \stackrel d\rightarrow C^{-1} \stackrel d\rightarrow C^{0}\stackrel d\rightarrow C^{1} \stackrel d\rightarrow\cdots $$
is a cochain complex in $\Pol_d$ then one may find an injective cochain map 
$C^\bullet\hookrightarrow J^\bullet$
with each $J^i$ injective and $J^i$ zero when $C^i$ is zero.
This is clear when $C^\bullet$ is an easy complex like $\cdots\to0\to F\to0\to\cdots$ or 
$\cdots\to0\to F\stackrel\id\to F\to 0\cdots$. Any $C^\bullet$ can be embedded into a direct sum of such easy complexes.

Recall that a cochain map $f:C^\bullet\to D^\bullet$ is called a \emph{quasi-isomorphism} if each $H^i(f):H^i(C^\bullet)\to
H^i(D^\bullet)$ is an isomorphism.

If $C^\bullet$ is a cochain complex in $\Pol_d$ that is \emph{bounded below}, meaning that $C^j=0 $ for $j\ll0$,
 then one may find a quasi-isomorphism $C^\bullet\hookrightarrow J^\bullet$
with each $J^j$ injective and $J^j$ zero when $j\ll0$. One may construct $J^\bullet$ as the total complex of a double
complex $K^\bullet_\bullet$ obtained by coresolving  like this:
Construct an exact complex of complexes $0\to C^\bullet\to K_0^\bullet\to K_1^\bullet\to\cdots$
where the 
$K_i^\bullet$ are complexes of injectives with $K_i^j=0$ when $C^j=0$. (Our double complexes commute so that a total complex
requires appropriate signs.)
One calls $C^\bullet\hookrightarrow J^\bullet$, or simply $J^\bullet$, an \emph{injective coresolution} of $C^\bullet$.  
Notice that we prefer our injective coresolutions to be injective as cochain maps, as indicated by the symbol $\hookrightarrow$. 
But \emph{any}\/ quasi-isomorphism $C^\bullet\rightarrow J^\bullet$  with each $J^j$ injective  is  called an injective
coresolution of $C^\bullet$.

If $f:J^\bullet\to \tilde J^\bullet $ is a quasi-isomorphism of bounded
below complexes of injectives, then the mapping cone of $f$ is a bounded below acyclic complex of injectives, hence split
and contractible,
and $f$ is a homotopy equivalence \cite[Proposition 0.3, Proposition 0.7]{Brown}.

\begin{Remark}
 Actually $\Pol_d$ has \emph{finite global dimension} \cite{Weibel} by \cite[A.11]{Jantzen book} 
so that even for an unbounded complex $C^\bullet$ there is 
a quasi-isomorphism $C^\bullet\hookrightarrow J^\bullet$
with each $J^j$ injective. Indeed the coresolution $0\to C^\bullet\to K_0^\bullet\to K_1^\bullet\to\cdots$ may be terminated
and thus one may use the total complex of the finite width double complex $K_0^\bullet\to K_1^\bullet\to\cdots\to K_M^\bullet$.
Also, a bounded complex is quasi-isomorphic to a bounded complex of injectives. (A complex $C^\bullet$ is called \emph{bounded}
if $C^i$ vanishes for $|i|\gg 0$.)
Passing to Kuhn duals one also sees that a bounded complex is quasi-isomorphic to a complex of 
projectives that is bounded.
\end{Remark}

\begin{Exercise}\label{quasi ex}
Let $C^\bullet \hookrightarrow D^\bullet $ be a quasi-isomorphism and let $C^\bullet\to E^\bullet$ be a cochain map.
Then $E^\bullet\hookrightarrow (D^\bullet\oplus E^\bullet)/C^\bullet $ is a quasi-isomorphism. There is a commutative diagram 
$$\begin{array}{ccc}
C^\bullet&\hookrightarrow&D^\bullet\\
\downarrow&&\downarrow\\
E^\bullet&\hookrightarrow&F^\bullet
\end{array}$$
with $E^\bullet\hookrightarrow F^\bullet$ a quasi-isomorphism. (Hint: One may also construct an anticommutative square.)
We will refer to this diagram as the \emph{base change} diagram.

Let $\KP_d$ be the homotopy category \cite[Exercise 1.4.5]{Weibel}
of bounded below cochain complexes   in 
$\Pol_d$. If $f:J^\bullet\to C^\bullet $ is a quasi-isomorphism of bounded
below complexes and the $J^i$ are injectives, then $f$ defines a split monomorphism in $\KP_d$.

The injective coresolution of a bounded below complex is unique up to homotopy equivalence. Here one does not require
the coresolutions to be injective as cochain maps. (But to prove it, consider a pair of injective coresolutions
with at least one of the two
injective as cochain map. Then use base change and coresolve.) 

Let $A^\bullet$ be an exact bounded
below complex. Its injective coresolutions are contractible. If $f:A^\bullet\to \tilde J^\bullet$ is a cochain map and
$\tilde J^\bullet $ is a bounded
below complexes of injectives, then $f$ is homotopic to zero. 
(Hint: Take $C^\bullet=A^\bullet$ and $E^\bullet=\tilde J^\bullet$ in the base change diagram.)

Let $f:J^\bullet\to \tilde J^\bullet $ be a morphism of bounded
below complexes of injectives. 
If there is a quasi-isomorphism $g:C^\bullet \hookrightarrow J^\bullet $ so that $fg=0$, then 
$f$ factors through $J^\bullet /C^\bullet$ and is thus homotopic to zero.
\end{Exercise}

\begin{Definition}\label{quasi zigzags}
We say that two complexes $C^\bullet$, $D^\bullet$ are \emph{quasi-isomorphic} if there are complexes 
$E^\bullet_0,\dots, E^\bullet_{2n}$  
and quasi-isomorphisms 
$f_i:E^\bullet_{2i}\to E^\bullet_{2i+1}$, $g_i:E^\bullet_{2i+2}\to E^\bullet_{2i+1}$
with $E^\bullet_0=C^\bullet$,  $E^\bullet_{2n}=D^\bullet$.
Thus $C^\bullet$, $D^\bullet$ are joined by zigzags of quasi-isomorphisms 
$E^\bullet_{2i}\to E^\bullet_{2i+1}\leftarrow E^\bullet_{2i+2}$.
\end{Definition}

It follows from exercise \ref{quasi ex} that injective coresolutions of quasi-isomorphic bounded below 
complexes are homotopy equivalent. This fact will underlie our choice of model for the derived category 
({\it cf.} \cite[Theorem 10.4.8]{Weibel}).

\section{Precomposition by Frobenius}
We will need the derived category $\DbP d$ to discuss the formality conjecture of 
Cha\l upnik, which is formulated in terms of $\DbP d$. 
In \ref{bif collapse} we will turn to the  collapsing conjecture of Touz\'e. Its formulation and proof do
not need anything about derived categories. That part of the story can be told entirely on the level of spectral
sequences of bicomplexes, but we leave it to the reader to disentangle the derived categories from the spectral sequences.
We find  the analogy between the 
formality problem and the collapsing conjecture instructive. 
Our use of derived categories is rather basic. We will model the derived category $\DbP d$ by a certain homotopy category
$\KbI_d$. We
could have phrased almost everything in terms of that homotopy category, but we like derived
categories and their intimate connections with spectral sequences.

On closer inspection the reader will find that even the existence of $\DbP d$ is not essential for the heart of the
arguments. One may simply view $\DP d$ as a source of inspiration and notation. 

\subsection{Derived categories}
One gets the \emph{derived category} $\DP d$ from the category of bounded below cochain complexes in $\Pol_d$ by forcing 
quasi-isomorphic complexes to be isomorphic. There are several ways to do that.
The usual way  is by formally inverting the 
quasi-isomorphisms. (The objects of the category do not change. Only the morphism sets are changed when throwing in
formal inverses.)

In our case there is a good alternative:
Replace every complex by an injective coresolution, then compute up to homotopy.
In fact, if $C^\bullet$, $D^\bullet$ are bounded below cochain complexes and $D^\bullet\hookrightarrow\tilde J^\bullet$ is
an injective coresolution, then $\Hom_{\DP d}(C^\bullet,D^\bullet)$ may be identified by \cite[Cor 10.4.7]{Weibel} with 
$\Hom_{\KP_d}(C^\bullet,\tilde J^\bullet)$,
where $\KP_d$ is the homotopy category \cite[Exercise 1.4.5]{Weibel}
of bounded below cochain complexes   in 
$\Pol_d$. One does not need to coresolve $C^\bullet$ here, but one may coresolve it too. If $J^\bullet$ is an injective 
coresolution of $C^\bullet$, then $\Hom_{\KP_d}(C^\bullet,\tilde J^\bullet)$ is isomorphic to
$\Hom_{\KP_d}(J^\bullet,\tilde J^\bullet)$.

Note that a cochain map
$f:C^\bullet\to D^\bullet$ is homotopic to zero if and only if it factors through the mapping cone of 
$\id:D^\bullet\to D^\bullet$. And this mapping cone is quasi-isomorphic to the zero complex.
Taking into account the $k$-linear structure it is thus not surprising that inverting quasi-isomorphisms
forces homotopic cochain maps to
 become equal
 \cite[Examples 10.3.2]{Weibel}.
The derived 
category may also be described \cite[10.3]{Weibel} by first passing to $\KP_d$ and then inverting the 
quasi-isomorphisms. 

Consider the full subcategory $\KI_d$ of the homotopy category whose objects are bounded below complexes of
injectives in $\Pol_d$. It maps into $\DP d$ and $\KI_d\to\DP d$ 
is an equivalence of categories \cite[Theorem 10.4.8]{Weibel}. 
One retracts $\DP d$ back to $\KI_d$ by sending a complex $C^\bullet$ to
an injective coresolution $J^\bullet$ of $C^\bullet$.
We use $\KI_d$ as our  working definition  of $\DP d$. Note that the definition of $\KI_d$ is
easy. No formal inverting is needed.
The way it works is that, when we try to understand morphisms in $\DP d$, we may model  an object of $\DP d$ 
by means of its image under the retract.

We view $\Pol_d$ as a subcategory of
$\DP d$ in the usual way: Associate to $F\in\Pol_d$ [an injective coresolution of]
the complex $\cdots\to0\to F\to0\to\cdots$ with $F$ in degree zero.
Write the complex as $F[-m]$ when $F$ is placed in degree $m$ instead.
The derived category encodes $\Ext$ groups as follows \cite[10.7]{Weibel}. If $F$, $G\in \Pol_d$ then 
$$\Hom_{\DP d}(F[-m],G[-n])\cong \Ext_{\Pol_d}^{m-n}(F,G).$$

One also has the \emph{bounded derived category}
$\DbP d$ which we think of as the full subcategory of $\DP d$ whose objects $C$ have vanishing $H^i(C)$ for $i\gg0$.
Let $\KbI_d$ be the subcategory of $\KI_d$ whose objects are homotopy equivalent to a bounded complex of
injectives in $\Pol_d$. Then $\KbI_d$ is our working definition of $\DbP d$.
(Recall that $\Pol_d$ has finite cohomological dimension.)

\subsection{The adjoint of the twist}
We now aim for a formality conjecture of Cha\l upnik \cite{Chalupnik}
 related to the collapsing conjecture of Touz\'e
 \cite[Conjecture 8.1]{Touze 2012}. 
These conjectures imply a powerful formula (Exercise \ref{expansion by formality}) for the
effect of Frobenius twist on Ext groups in the category of strict polynomial functors. 
For the application to the (CFG) conjecture 
we will need to extend the theory from strict polynomial functors to strict polynomial bifunctors, 
but the difficulties are already
visible for strict polynomial functors.

Let $A\in \Pol_e$. The example we have in mind is $A=I^{(r)}$, the $r$-th Frobenius twist.
Precomposition with $A$ defines a functor $\Pol_d\to \Pol_{de}: F\mapsto F\circ A$.
So the example we have in mind is $F\mapsto F^{(r)}$. The functor $F\mapsto F\circ A$ extends to a functor 
${-}\circ A:\KbI_d\to \DbP {de}$, hence a functor $\DbP d\to \DbP {de}$. We first seek its right adjoint $\KrA$.
For an object $J^\bullet$ of $\KbI_{de}$ put 

$$\KrA(J^\bullet)(V):=\Hom_{\Pol_{de}}(\Gamma^{dV}\circ A,J^\bullet),$$ 
where the right hand side is
viewed as a complex in $\Pol_{d}$ of functors $V\mapsto \Hom_{\Pol_{de}}(\Gamma^{dV}\circ A,J^i)$.
Observe that this complex is homotopy equivalent to a bounded complex.
If $G\in \DbP {de}$, then we take an injective coresolution $J^\bullet$ of $G$ and put 
$\KrA(G):=\KrA(J^\bullet)$.
Our claim is that $$\Hom_{\DbP {de}}(F\circ A,G)=\Hom_{\DbP {d}}(F,\KrA(G)),$$
for $F\in \DbP {d}$, $G\in \DbP {de}$.

Now take $Z\in \GVde$ of dimension at least $de$. Then every object in $\DP {de}$ is quasi-isomorphic to one of the form
$$G=\cdots\to k^{n_i}\otimes_kS^{de}_Z\to k^{n_{i+1}}\otimes_kS^{de}_Z\to \cdots,$$
so we may assume that $G$ is actually of this form. Notice that $k^{n_i}\otimes_kS^{de}_Z\to k^{n_{i+1}}\otimes_kS^{de}_Z$
is given by an $n_{i+1}$ by $n_i$ matrix with entries in $\End_{\GVde}(Z)$.
We may also assume $F=F^\bullet$ consists of projectives and is bounded. 
(We wish to use \emph{balancing} \cite[2.7]{Weibel},
which is the principle that both projective resolutions and injective coresolutions may be used to compute `hyper Ext'.
See also section \ref{untwisting collapse}.
We do not use an injective coresolution of $F$.)
Put $F_i=F^{-i}$.
Now $\Hom_{\DP {de}}(F\circ A,G)$ is computed as the $H^0$ of the total complex associated to the bicomplex 
$\Hom_{\Pol_{de}}(F_i\circ A,G^j)$ and $\Hom_{\DbP {d}}(F,\KrA(G))$ is similarly computed by way of a 
bicomplex \cite[2.7.5, Cor 10.4.7]{Weibel}.
So let us compare the bicomplexes. We have 
$\Hom_{\Pol_{de}}(F_i\circ A,G^j)=\Hom_{\Pol_{de}}(F_i\circ A,k^{n_j}\otimes_kS^{de}_Z)=
k^{n_j}\otimes_k\Hom_{\Pol_{de}}(F_i\circ A,S^{de}_Z)=k^{n_j}\otimes_k(F_i\circ A)^\#Z$ and 
$\Hom_{\Pol_{d}}(F_i,\KrA(G)^j)=\Hom_{\Pol_{d}}(F_i,
V\mapsto \Hom_{\Pol_{de}}(\Gamma^{dV}\circ A,k^{n_j}\otimes_kS^{de}_Z))=
k^{n_j}\otimes_k\Hom_{\Pol_{d}}(F_i,(V\mapsto (\Gamma^{dV}\circ A)^\#Z))=
k^{n_j}\otimes_k\Hom_{\Pol_{d}}(F_i,S^d_{A^\#Z})=k^{n_j}\otimes_kF_{i}^\# A^\#Z$.
The claim follows. (Exercise.)

\begin{Remark} These bicomplexes $\Hom_{\Pol_{de}}(F_i\circ A,G^j)$ and $\Hom_{\Pol_{d}}(F_i,\KrA(G)^j)$
are meaningful by themselves. The fact that their total complexes are quasi-isomorphic may also serve as
motivation for the definition of $\KrA(G)$. This does not explicitly involve the derived category.
It is closer in spirit to section \ref{untwisting collapse}. The bicomplexes do not require that $F^\bullet$
is a bounded complex of projectives. A bounded above complex of projectives would do.
\end{Remark}

\subsection{Formality}\label{formality section} A bounded below cochain complex $C^\bullet$ in $\Pol_d$
is called \emph{formal} if it is isomorphic in the derived category $\DP d$
to a complex $E^\bullet$ with zero differential. 
Notice that here we do not replace $E^\bullet$ with an injective coresolution,
because that usually spoils the vanishing of the differential. 
Notice also that $E^i\cong H^i(E^\bullet)\cong H^i(C^\bullet)$.
One can show that the isomorphism in the derived category is given by
a single zigzag of quasi-isomorphisms $C^\bullet\to D^\bullet \leftarrow E^\bullet$ or, dually, 
a single zigzag of quasi-isomorphisms $C^\bullet\leftarrow D^\bullet \to E^\bullet$.
One may also define a cochain complex to be formal if it is quasi-isomorphic in the sense of definition
\ref{quasi zigzags} to a complex with zero differential. So one does not need the derived category to
introduce formality.

\begin{Exercise}A complex has differential zero if
and only if it is a direct sum of complexes each of which is concentrated in one degree.

Let $m$ be an integer.
Let $C^\bullet$ be a bounded below cochain  complex with $H^i(C^\bullet)=0$ for $i\neq m$.
Show that $C^\bullet$ is formal by constructing a zigzag of quasi-isomorphisms 
$C^\bullet\hookleftarrow D^\bullet \twoheadrightarrow E^\bullet$
where $D^i=0$ for $i>m$, $E^i=0$ for $i\neq m$.
\end{Exercise}

\begin{Remark}\label{not formal}
 Let the $2$-fold extension 
$0\to F\to G \stackrel f\to H \to K\to 0$ represent \cite[2.6]{Benson I}
a nonzero element of $\Ext_{\Pol_d}^2(K,F)$. 
One can show that $$\cdots\to 0\to G \stackrel f\to H \to 0\to\cdots$$
is not formal. 

%

For example,
consider the $2$-fold extension
$$0\to I^{(1)}\to S^p\stackrel\alpha\to \Gamma^p\to I^{(1)}\to 0$$
of \cite[Lemma 4.12]{Friedlander-Suslin}  where $\alpha_V:S^p(V)\to \Gamma^p(V)$ is the symmetrization 
homomorphism, $\alpha_V(v_1\cdots v_p)=\sum_{\sigma\in\cS_p}v_{\sigma^{-1}(1)}\otimes\cdots\otimes v_{\sigma^{-1}(p)}$.
It represents a nontrivial class by \cite[Lemma 4.12, Theorem 1.2]{Friedlander-Suslin}. So 
$$\cdots\to0\to  S^p\stackrel\alpha\to \Gamma^p\to 0\to \cdots$$ is not formal. See also Exercise \ref{proof not formal}.
\end{Remark}

Now let $A=I^{(r)}$. Then we write $\KrA$ as $\Kr$. Let $E_r$ be the graded vector space of dimension $p^r$ which equals
$k$ in dimensions $2i$, $0\leq i<p^r$. We view any graded vector space also as a cochain complex with zero differential 
and as a 
$\Gm$-module with weight $j$ in degree $j$. For example, if $G \in \Pol_d$ then the $\Gm$ action on $E_r$ induces
one on $G_{E_r}$ so that $G_{E_r}$
is graded and thus a
complex with differential zero.

A conjecture of  Cha\l upnik, now says
\begin{Conjecture}[Formality]\label{formality conjecture}
For $G$ in $\Pol_d$ one has $\Kr(G^{(r)})\cong G_{E_r}$ in $\DbP d$. In particular,  $\Kr(G^{(r)})$ is formal.
\end{Conjecture}
This is a variant of the collapsing conjecture of Touz\'e \cite[Conjecture 8.1]{Touze 2012}. It is stronger.
Both conjectures are theorems now \cite{Chalupnik}, \cite{Touze 2013}, \ref{untwisting formality}.

Note that $G_{E_r}$ is formal by definition. Note also that all its weights are even. 
So as a cochain complex it lives in even degrees. That already implies formality.

\begin{Exercise}\label{expansion by formality}(Compare \cite[Corollary 5]{Touze 2013}.)
 Let $F$, $G\in\Pol_d$. Assuming the formality conjecture, define a grading on $\Ext^\bullet(F,G_{E_r})$ so that its degree $i$ 
subspace is isomorphic to  
$\Ext_{\Pol_{p^rd}}^i(F^{(r)},G^{(r)})$.
\end{Exercise}

\begin{Remark}Let $J^\bullet$ be a 
bounded injective coresolution of $ G^{(r)}$.
 If one can show formality of $\Kr(J^\bullet)$, then one can also show it is quasi-isomorphic to $G_{E_r}$.
The main problem is formality of $\Kr(J^\bullet)$.  This problem does not require derived categories. 
But the problem is easier to motivate
in the language of derived categories.
\end{Remark}

Following suggestions by Touz\'e let us give some evidence for the formality conjecture in the simplest case: $p=2$, $r=1$. 
Instead of $\Kr(G^{(1)})$ we will study $\Kr(G^{(1)})\circ I^{(1)}$ and show that it is formal. 
So we will be off by one Frobenius twist. While we know how to untwist in $\Pol_d$ context, something  more 
will be needed to do  untwisting
in $\DbP d$ context. We postpone this issue until \ref{untwisting formality}. 

Now $\Kr(G^{(1)})\circ I^{(1)}$ is represented by the complex
$\Hom_{\Pol_{2d}}(\Gamma^{dV^{(1)}}\circ I^{(1)},J^\bullet)$ in $\Pol_{2d}$,
 where $J^\bullet$ is a bounded injective coresolution of
$G^{(1)}$. Observe that $\Gamma^{dV^{(1)}}\circ I^{(1)}=(\Gamma^{d{(1)}})^V$.
This is where the extra twist helps: It turns out that $(\Gamma^{d{(1)}})^V$ is easier than $\Gamma^{dV}\circ I^{(1)}$.
Rewrite our complex as $\Hom_{\Pol_{2d}}(J^{\bullet\#},(S^{d{(1)}})_V)$. 
We first recall a standard injective coresolution of $(S^{d{(1)}})_V$.

\subsection{A standard coresolution in characteristic two}
It is here that the assumptions on $p$ and $r$ help. In general one needs the Troesch complexes to see that
$\Kr(G^{(r)})\circ I^{(r)}$ equals $G_{E_r}\circ I^{(r)}$ in $\DbP {p^rd}$
and we refer to 
\cite{Touze 2012}, \cite{Touze 2013} for details. 

Let $T$ be the group scheme of diagonal matrices in $\GL_2$.
If $W\in \V_k$ then $T$ acts through $k^2$ on the symmetric algebra $S^*(k^2\otimes_k W)$ with weight space 
$S^i(W)\otimes_k S^{j}(W)$ of
weight $(i,j)$. So the $S^i(W)\otimes_k S^{j}(W)$ are direct summands of $S^{i+j}_{k^2}(W)$ and $S^i\otimes_kS^j$ is an injective in
$\Pol_{i+j}$ because it is a summand of an injective. Now recall $p=2$.
We make the algebra $S^*(k^2\otimes_k W)=S^*W\otimes_kS^*W$ into a differential graded algebra with 
differential $d$ whose restriction to
$S^1(k^2\otimes_k W)$ is given by $\pmatrix{0&0\cr 1&0}\otimes\id_W$.
So if $W$ has dimension one then the differential graded algebra is isomorphic to the polynomial ring $k[x,y]$
in two variables and the differential is $y\frac{\partial}{\partial x}$.
The subcomplex 
$S^{2n}W\to S^{2n-1}W\otimes_kS^1W\to\cdots\to S^{2n-i}W\otimes_kS^iW\to\cdots S^{1}W\otimes_k S^{2n-1}W\to S^{2n}W$ 
is a coresolution 
of $S^nW^{(1)}$. (Exercise. Use the exponential property.)
So we have a standard coresolution 
$$ S^{2n}\to S^{2n-1}\otimes_kS^1\to\cdots\to S^{2n-i}
\otimes_kS^i\to\cdots S^{1}\otimes S^{2n-1}\to S^{2n}\to0$$ of $S^{n(1)}$.
We now coresolve $(S^{d{(1)}})_V$ by 
$$R_V^{2d\bullet}:\quad S_V^{2d}\to S_V^{2d-1}\otimes_k S_V^{1}\to\cdots S_V^{2d-i}\otimes_k S_V^{i}\to
\cdots\to S_V^{1}\otimes_k S_V^{2d-1}\to S_V^{2d}.$$ 

\begin{Exercise}\label{proof not formal}
Keep $p=2$. Determine $\Hom_{\Pol_2}(S^2,S^2)$ and $\Hom_{\Pol_2}(\Gamma^{2},I\otimes I)$.
 Show there is no nonzero cochain map from the complex
$$C^\bullet:\cdots\to0\to S^2\stackrel\alpha\longrightarrow \Gamma^2\longrightarrow 0\to0\to\cdots$$ of Remark \ref{not formal} 
to the injective coresolution
$$D^\bullet:\cdots\to0\to S^2\to I\otimes I\to S^2\to0\to\cdots$$ of $I^{(1)}$. 
(In both complexes the first $S^2$ is placed in degree 0.)
By \cite[Cor 10.4.7]{Weibel} this implies there is no nonzero morphism from $C^\bullet$ to $D^\bullet$ in the derived category
$\DP 2$.
Confirm the claim in Remark \ref{not formal} that $C^\bullet$ 
is not formal.
\end{Exercise}

\subsection{Formality continued}
We are studying the complex $\Hom_{\Pol_{2d}}(J^{\bullet\#},(S^{d{(1)}})_V)$ in $\Pol_{2d}$ up to quasi-isomorphism.
We may replace it with the total complex of the double complex $\Hom_{\Pol_{2d}}(J^{\bullet\#},R_V^{2d\bullet})$ and then 
(by `balance' \cite[2.7]{Weibel}) with the 
complex $\Hom_{\Pol_{2d}}(G^{{(1)}\#},R_V^{2d\bullet})$ in $\Pol_{2d}$.
If one forgets the differential then this is just 
$\Hom_{\Pol_{2d}}(G^{{(1)}\#},S_{k^2\otimes_kV}^{2d})=G^{{(1)}}(k^2\otimes_kV)$ and we now inspect its weight spaces 
for our 
torus $T$. Because of the Frobenius twist in $G^{(1)}$ the weights are all multiples of $p$, and $p$ equals $2$ now.
On the other hand, on $\Hom_{\Pol_{2d}}(G^{{(1)}\#},S_V^{2d-i}\otimes_k S_V^{i})$ the weight is simply $(2d-i,i)$.
So the only nonzero terms in the complex $\Hom_{\Pol_{2d}}(G^{{(1)}\#},R_V^{2d\bullet})$ are in even degrees and formality 
follows.
Moreover, in even degree $2i$ one gets the weight space of degree $(2d-2i,2i)$ of  $G^{{(1)}}(k^2\otimes_kV)$.
Let $\Gm$ act on $k^2$ with weight zero on $(1,0)$ and weight one on $(0,1)$.
So now $E_1=(k^2)^{{(1)}}$ as $\Gm$-modules.
As a (graded) functor in $V$ we get that $\Hom_{\Pol_{2d}}(G^{{(1)}\#},R_V^{2d\bullet})$ is $(G^{{(1)}})_{k^2}$ or 
$G_{E_1}\circ I^{(1)}$. So we have seen that for $p=2$, $r=1$ the complex $\Kr(J^\bullet)\circ I^{(1)}$ is quasi-isomorphic to
$G_{E_r}\circ I^{(1)}$. That means that $\Kr(G^{(r)})\circ I^{(1)}\cong G_{E_r}\circ I^{(1)}$ in $\DP {2d}$ for $p=2$, $r=1$. 

Now we would like to untwist to get  the formality conjecture for $p=2$, $r=1$.
It is not obvious how to do that. One needs constructions with better control of the functorial 
behavior. In his solution in \cite{Touze 2010} of the (CFG) conjecture
Touz\'e faced similar difficulties.
As standard coresolutions he used Troesch coresolutions. 
They are not functorial. This is the main obstacle that he had to get around in order to construct the classes $c[m]$. 
His approach in \cite{Touze 2010} is to invent a new category, the twist-compatible category,
on which the Troesch construction is functorial and which is just big enough to contain a repeated reduced bar construction
that coresolves divided powers. 

In the proof \cite{Touze 2013} of the collapsing conjecture a different argument is used. We call it untwisting the collapse of 
a hyper Ext spectral sequence. It comes next.

\subsection{Untwisting the collapse of a hyper Ext spectral sequence}\label{untwisting collapse}
Let $C^\bullet$ be a bounded above complex in $\Pol_d$ and $D^\bullet$ a bounded below complex in
$\Pol_d$. Put $C_i=C^{-i}$.
Let $J^\bullet$ be a bounded below complex of injectives that coresolves $D^\bullet$.
The homology groups of the total complex
$\Tot\Hom_{\Pol_d}(C_{\bullet},J^\bullet)$ of the bicomplex $\Hom_{\Pol_d}(C_{i},J^j)$ are known as 
\emph{hyper Ext groups} of $C^\bullet$, $D^\bullet$. Consider the second 
 spectral sequence of the bicomplex
$\Hom_{\Pol_d}(C_{i},J^j)$
$$E_2^{ij}=H^i\Hom_{\Pol_d}(H_{j}(C_\bullet),J^\bullet)\Rightarrow H^{i+j}\Tot\Hom_{\Pol_d}(C_{\bullet},J^\bullet).$$
We call it 
the  \emph{hyper Ext spectral sequence} associated with $(C^\bullet, D^\bullet)$. 
It is covariantly functorial in $D^\bullet$ and contravariantly functorial   in $C^\bullet$. 
Say $\tilde C^\bullet\to C^\bullet$, $D^\bullet\to\tilde D^\bullet$
are quasi-isomorphisms. Let $\tilde J^\bullet$ be an injective coresolution of $\tilde D^\bullet$.
Then $J^\bullet$, $\tilde J^\bullet$ are quasi-isomorphic complexes of injectives, hence homotopy equivalent.
The 
hyper Ext spectral sequence associated with $(C^\bullet, D^\bullet)$ is isomorphic with the
 hyper Ext spectral sequence associated with $(\tilde C^\bullet, \tilde D^\bullet)$. (Check this.)
In particular, if $C^\bullet$ is formal, then the spectral sequence is a direct sum of spectral sequences with
just one row, so that the spectral sequence degenerates at page two. We also say that the spectral sequence \emph{collapses}.

Now suppose that we do not know that $C^\bullet$ is formal, but only that $C^{\bullet(1)}$ is formal. 
Frobenius twist $G\mapsto G^{(1)}$ defines an embedding of $\Pol_d$ into $\Pol_{dp}$. 
Coresolving $J^{\bullet(1)}$ we get a map from the hyper Ext
spectral sequence $E$ of $(C^\bullet,D^\bullet)$ to the hyper Ext
spectral sequence $\tilde E$ of $(C^{\bullet(1)}, D^{\bullet(1)})$.
Now we make the extra assumption that $D^\bullet$ is concentrated in one degree.
Say degree zero, to keep notations simple. 
Write $D^\bullet$ as $D$. Then the second page of $E$ is given by $E_2^{ij}=\Ext_{\Pol_d}^i(H_{j}(C_\bullet),D)$
and the second page of $\tilde E$ is given by $\tilde E_2^{ij}=\Ext_{\Pol_{dp}}^i(H_{j}(C_\bullet)^{(1)},D^{(1)})$.
Now  the map $ E_2^{ij}\to \tilde E_2^{ij}$ is injective by the Twist Injectivity Theorem \ref{twist injectivity}.
We conclude that $E$ itself degenerates at page two by means of the following basic lemma about spectral sequences.

\begin{Lemma}
Let $E\to \tilde E$ be a morphism of spectral sequences that is injective at the second page. If $\tilde E$
degenerates at page two, then so does $E$.
\end{Lemma}

\subsubsection*{Proof}
The second page $E_2$ of $E$ with differential $d_2$ may be viewed as a subcomplex of $\tilde E_2$ with differential
$\tilde d_2$.
So the differential $d_2$ of $E_2$ vanishes and $E_3\to \tilde E_3$ is also injective. But the differential of
 $\tilde E_3$ vanishes, so the differential of $E_3$ vanishes again. Repeat. \qed

\

So we do \emph{not} need the formality of  $C^{\bullet}$ to conclude the collapsing of 
$E_2^{ij}=\Ext_{\Pol_d}^i(H_{j}(C_\bullet),D)\Rightarrow H^{i+j}\Tot\Hom_{\Pol_d}(C_{\bullet},J^\bullet)$.
Formality of $C^{\bullet(1)}$ suffices.
We have `untwisted' the collapsing.

\begin{Exercise}(Untwisting Formality)\label{untwisting formality}
Let $C^\bullet$ be a bounded  complex in $\Pol_d$ such that $C^{\bullet(1)}$ is formal and let 
$D^\bullet$ a formal bounded below complex in
$\Pol_d$. Show that the hyper Ext
spectral sequence $E$ of $(C^\bullet,D^\bullet)$ collapses.
Say $C^\bullet$ has nonzero cohomology.
Put $m=\min\{\;i\mid H^i(C^\bullet)\neq0\;\}$.
Now take for $D^\bullet$ the maximal subcomplex of $C^\bullet$ with $D^i=0$ for $i>m$.
Then $H^m(D^\bullet)\to H^m(C^\bullet)$ is an isomorphism
and $H^i(D^\bullet)=0$ for $i\neq m$. The complex $D^\bullet$ is formal. 
Let $D^\bullet\to J^\bullet$ be an injective coresolution again. 
Recall that $H^0(\Tot\Hom_{\Pol_d}(C_{\bullet},J^\bullet))=\Hom_{\KP_d}(C^\bullet,J^\bullet)$ by \cite[2.7.5]{Weibel}.
Compare the collapsed hyper Ext spectral sequences of  $(C^\bullet,D^\bullet)$ and $(D^\bullet,D^\bullet)$.
Show that $\Hom_{\KP_d}(C^\bullet,J^\bullet)\to \Hom_{\KP_d}(D^\bullet,J^\bullet)$ is surjective.
Choose $f:C^\bullet\to J^\bullet$ so that the composite $D^\bullet\to C^\bullet\stackrel{ f}{\to} J^\bullet$ is a 
quasi-isomorphism.
Show that $C^\bullet$ is quasi-isomorphic to $J^\bullet\oplus (C^\bullet/D^\bullet)$. 
Show that $C^\bullet$ is formal by induction on the number of nontrivial $H^i(C^\bullet)$.
This establishes untwisting of formality. Notice that we did not mention $\DbP {d}$ in this exercise. But recall
that $\Hom_{\DP d}(C^\bullet,D^\bullet)$ may be identified with 
$\Hom_{\KP_d}(C^\bullet, J^\bullet)$.
\end{Exercise}

This finishes the proof of the formality conjecture for the case: $p=2$, $r=1$. 

\section{Bifunctors and CFG}
There are some more ingredients entering into the proof of the CFG theorem in \cite{TvdK}.
The paper \cite{TvdK} has an extensive introduction, which we recommend to the reader.
We now provide a companion to that introduction.

The proof of the CFG conjecture takes several steps.
First one reduces to the case of $\GL_n$. This uses a transfer principle, reminiscent of Shapiro's Lemma, 
that can be traced back to the nineteenth century.
Next one needs to know about Grosshans graded algebras and good filtrations.
The case were the coefficient algebra $A$ is a Grosshans graded algebra lies in between the general case and the case
of good filtration. In the good filtration case  
CFG is known by invariant theory. There is a spectral sequence connecting the Grosshans graded case
with the general case and another spectral sequence connecting it with the good filtration case.
We need to get these spectral sequences under control. That is done by finding an algebra of operators, operating on the
spectral sequences, and establishing finiteness properties of the spectral sequences  with respect to the operators.
It is here that the classes of Touz\'e come in. They allow a better grip on the operators. 

Now we introduce some of these notions.
\subsection{Costandard modules}
Let $k=\F_p$  and put $G=\GL_n$, $n\geq2$.
We have already introduced the torus $T$ of diagonal matrices. Our standard Borel group $B$ will be 
the subgroup scheme with $B(R)$ equal to the subgroup of upper triangular matrices of $G(R)$.
Similarly $U$, the unipotent radical of $B$, is the subgroup scheme with $U(R)$ equal to the 
subgroup of upper triangular matrices with ones on the diagonal.
The \emph{Grosshans height} $\hgt$, also known as the sum of the coroots associated to the positive roots, is
given by $$\hgt(\lambda)=\sum_{i<j}\lambda_i-\lambda_j=\sum_i(n-2i+1)\lambda_i.$$ 
Here we use the ancient convention that the roots of $B$ are positive.
If $V$ is a representation 
of $G$, let us say that it has \emph{highest weight} $\lambda$ if $\lambda$ is a weight of $V$ and all other weights 
$\mu$
have strictly
smaller Grosshans height $\hgt(\mu)$. (This nonstandard convention is good enough for the present purpose.)
Irreducible $G$-modules have a highest weight and are classified up to isomorphism by that weight.
Write $L(\lambda)$ for the irreducible module with highest weight $\lambda$. The weight space of weight $\lambda$
in $L(\lambda)$ is one dimensional and equal to the subspace $L(\lambda)^U$ of $U$-invariants.

We now switch to geometric language as if we are dealing with varieties. 
In other words, we switch from the setting of group schemes \cite{Demazure-Gabriel}, \cite{Jantzen book},
\cite{Waterhouse book}
to algebraic groups and varieties defined over $\F_p$ \cite{Springer LAG}.

The \emph{flag variety} \cite[8.5]{Springer LAG} $G/B$ is a \emph{projective variety} \cite[I \S2]{Hartshorne book},
 \cite[1.7]{Springer LAG}, 
not an affine variety. 
Given $L(\lambda)$ as above
there is an \emph{equivariant line bundle} \cite[8.5.7]{Springer LAG} $\cL_\lambda$
on $G/B$ so that its module $\nabla(\lambda)$ of global \emph{sections} \cite[8.5.7-8]{Springer LAG} on 
$G/B$ has a unique irreducible submodule, and this submodule is isomorphic to
$L(\lambda)$. The \emph{costandard module} $\nabla(\lambda)$ is finite dimensional (because $G/B$ is a projective variety).
 The weight space of weight $\lambda$
in $\nabla(\lambda)$ is also one dimensional and equal to the subspace $\nabla(\lambda)^U$ of $U$-invariants.
Every other $G$-module $V$ whose weight space $V_\lambda$ of weight $\lambda$ is one dimensional and equal to $V^U$ embeds 
into $\nabla(\lambda)$.
\emph{Kempf vanishing} \cite[II Chapter 4]{Jantzen book} says that $\cL_\lambda$ has no higher sheaf cohomology on $G/B$. 
One derives from this \cite{CPSvdK}
that 
$H^{>0}(G,\nabla(\lambda))$ vanishes. All nontrivial cohomology of $G$-modules is due to the distinction between the irreducible
modules $L(\lambda)$
and the {costandard} modules $\nabla(\lambda)$. The dimensions of the weight spaces of $\nabla(\lambda)$ 
are given by the famous 
Weyl character formula 

$$\Char(\nabla(\lambda))=\frac{\sum_{w\in W}(-1)^{\ell(w)}e^{w(\lambda+\rho)}}
{e^\rho\prod_{\alpha>0}(1-e^\alpha)}.$$
We do not explain the precise meaning here but just observe that the formula is characteristic free. 
The dimensions of
the weight spaces of $\nabla(\lambda)$ are the same as in the irreducible $\GL_n(\CC)$-module with highest weight $\lambda$.
Determining the dimensions of
the weight spaces of $L(\lambda)$ is less easy in general, to put it mildly.

\begin{Example}
 Let $V=k^n$ be the defining representation of $\GL_n$ over $\F_p$. The symmetric powers $S^m(V)$ are costandard modules.
More specifically, $S^m(V)$ is $\nabla((m,0,\cdots,0))$.
When $m=p^r$  the irreducible submodule $L((m,0,\cdots,0))$ of $S^m(V)$ is spanned by the $v^{p^r}$.
\end{Example}

If $V$ is a nonzero $G$-submodule of $\nabla(\lambda)$ then it determines a map $\phi_V$
from the flag variety $G/B$ to the projective
space whose points are codimension one subspaces of $V$, or one dimensional subspaces of $V^\vee$. 
(To a point of $G/B$ one associates the codimension one subspace of
$V$ consisting of sections vanishing at the point. Then one takes the elements in the dual that vanish on the codimension
one subspace.)
The image of $G/B$ under $\phi_V$
is isomorphic to $G/\tilde P$, where 
$\tilde P$ is the scheme
theoretic stabilizer of the image  of the point $B$. 
Here `scheme theoretic' indicates that the functorial 
interpretation of group schemes is needed. The image of the point $B$ is the highest weight space of $V^\vee$.
 The group scheme $\tilde P$ need not be reduced \cite{Lauritzen}, \cite{Wenzel}, 
but the image of $\tilde P$ under a sufficiently high power $F^r$ of
the Frobenius homomorphism $F:G\to G$ is  the stabilizer $P$
of the highest weight space of $\nabla(\lambda)^\vee$. This $P$ is an ordinary \emph{parabolic subgroup} 
\cite[6.2]{Springer LAG} and thus reduced, meaning
that its coordinate ring is reduced. There is a graded algebra associated with the image of  $\phi_V$. 
This algebra $A_V$ is known as \emph{coordinate ring of the affine cone} over
the image of $\phi_V$. It is a graded $k$-algebra, generated as a $k$-algebra by its degree one part, which is $V$.
This is typical for closed subsets of a projective space: Such a subset does not have an ordinary coordinate ring like an affine 
variety would, but a graded coordinate ring \cite[II Corollary 5.16]{Hartshorne book}.

Similarly one has a graded algebra 
$$A_{\nabla(\lambda)}=\bigoplus_{m\geq0}\Gamma(G/B,\cL_\lambda^m)$$
associated with the image $G/P$ of $G/B$ in the projective space whose
points are codimension one subspaces of $\nabla(\lambda)$.
The algebra $A_V$  may be embedded into  $A_{\nabla(\lambda)}$.
 Mathieu observed \cite[3.4]{Mathieu G} that the two affine cones have the same rational points over fields and concluded 
from this
that for $r\gg0$
the smaller algebra
 contains all $f^{p^r}$ for $f$ in the 
larger algebra.  This is not always
 the same $r$ as in $F^r$ above. 

\subsection{Grosshans filtration}
The situation above generalizes.
If $V$ is a possibly infinite dimensional $G$-module we define its \emph{Grosshans filtration} to be the filtration
$V_{\leq-1}=0\subseteq V_{\leq0}\subseteq V_{\leq1}\subseteq V_{\leq2}\cdots$ where $V_{\leq i}$ is the largest 
$G$-submodule of $V$ all whose weights
$\mu$ satisfy $\hgt(\mu)\leq i$. The associated graded $\bigoplus_i V_{\leq i}/V_{\leq i-1}$ we call the \emph{Grosshans graded}
$\gr V$. It can naturally be embedded into a direct sum $\hull_\nabla(\gr V)$ of costandard modules in such a way that
no new $U$-invariants are introduced: $(\gr V)^U=(\hull_\nabla(\gr V))^U$.
We say that $V$ has \emph{good filtration} \cite[II 4.16 Remarks]{Jantzen book} if $\gr V$ itself is a 
direct sum of costandard modules, 
in which case
$\gr V=\hull_\nabla(\gr V)$ \cite[Theorem 16]{Grosshans contr}. As costandard modules have no higher 
$G$-cohomology, a
module with good filtration has vanishing higher $G$-cohomology.
One says that a module has \emph{finite good filtration dimension} if it has a finite coresolution by modules with good
filtration. Such a module has only finitely many nonzero $G$-cohomology groups. 

If $A\in\Rg_k$ is a $k$-algebra with $G$-action, so that the multiplication map $A\otimes_k A\to A$ is a $G$-module map, then
$\gr A$ and $\hull_\nabla(gr A)$ are also $k$-algebras with $G$-action. Moreover, if $A$ is of finite type, then so are 
$\gr A$ and $\hull_\nabla(gr A)$ by Grosshans \cite{Grosshans contr}. 
And then there is an $r$ so that $\gr A$  contains all $f^{p^r}$ for $f$ in the 
larger algebra $\hull_\nabla(gr A)$. All higher $G$-cohomology of
$A$ is due to the distinction between $\gr A$ and $\hull_\nabla(\gr A)$.
It is here that Frobenius twists and Frobenius kernels enter the picture. (In this subject area a \emph{Frobenius kernel}
refers to the finite group scheme which is the scheme theoretic kernel of an iterated Frobenius map $F^r:G\to G$.)
In general we have no grip on the size of the minimal
$r$ so that
 $\gr A$  contains all $f^{p^r}$. This is where the results get much more qualitative than those of Friedlander and 
Suslin.

\begin{Problem}
 Given your favorite $A$, estimate the $r$ such that $\gr A$  contains all $f^{p^r}$ for $f$ in the 
larger algebra $\hull_\nabla(\gr A)$. Such an estimate is desirable because one may give a bound on the Krull 
dimension of $H^\even(G,A)$ in terms of $r$, $n$ and 
$\dim A$ by inspecting the proof in \cite{TvdK}. 
\end{Problem}

\subsection{The classes of Touz\'e}
The \emph{adjoint representation} $\gl_n$ of $\GL_n$ is defined as the $k$-module of $n$ by $n$ matrices over $k$ with 
$\GL_n(R)$ acting by conjugation on the set $M_n(R)=\gl_n\otimes_kR$ of $n$ by $n$ matrices over $R$.
This is also known as the adjoint action on the Lie algebra. The adjoint representation is not a polynomial representation
as soon as $n\geq2$. We now have all the ingredients to state the theorem of Touz\'e on \emph{lifted classes} proved
using strict polynomial bifunctors.
The base ring is our field $k=\F_p$ and $n\geq2$.

\begin{Theorem}[Touz\'e \cite{Touze 2010}. Lifted universal cohomology classes]\label{divided powers}
There are cohomology classes $c[m]$
so that
\begin{enumerate}
\item $c[1]\in H^{2}(\GL_{n},\gl_n^{(1)})$ is nonzero,
\item For $m\geq1$ the class $c[m]\in H^{2m}(\GL_{n},\Gamma^{m}(\gl_n^{(1)}))$ lifts 
$c[1]\cup\cdots\cup c[1]\in H^{2m}(\GL_{n},\bigotimes^{m}(\gl_n^{(1)}))$.
\end{enumerate}
\end{Theorem}

\subsection{Strict polynomial bifunctors}
The representations $\Gamma^{m}(\gl_n^{(1)}))$ in the `lifted classes' theorem of Touz\'e
are not polynomial. To capture their behavior one needs the \emph{strict polynomial bifunctors} of Franjou and Friedlander
\cite{Franjou-Friedlander}.
We already encountered them in disguise when discussing parametrized functors.
An example of a strict polynomial bifunctor is the bifunctor 
$$\Hom_{\GVd_k}({-_1},{-_2}):\GVd^\opp_k\times \GVd_k\to \V_k;\ (V,W)\mapsto \Hom_{\GVd_k}(V,W).$$
It is contravariant in ${-_1}$ 
 and covariant in ${-_2}$. 
More generally one could consider the category $\Pol^d_e$ of $k$-bilinear functors $\GVd^\opp_k\times \Gamma^e\V_k\to \V_k$.
Do not get confused by the strange notation $\Pol^\opp\times \Pol$ for $\bigoplus_{d,e}\Pol^d_e$ used in 
\cite{Franjou-Friedlander}. It is not a product.

If $\cal A$ and $\cal B$ are $k$-linear categories, then one can form the $k$-linear category ${\cal A}\otimes_k\cal B$ whose
objects are pairs $(A,B)$ with $A\in\cal A$, $B\in \cal B$.
For morphisms one puts $\Hom_{{\cal A}\otimes_k\cal B}((A,B),(A',B'))=\Hom_{{\cal A}}(A,A')\otimes_k\Hom_{\cal B}(B,B')$.
One may then define the category $\Pol^d_e$ of strict polynomial bifunctors of bidegree $(d,e)$ to be the category of $k$-linear
functors from $\GVd^\opp_k\otimes_k \Gamma^e\V_k$ to $\V_k$.

One gets more bifunctors by composition.
For instance, $\Gamma^{m}(\gl^{(1)})$ is the strict polynomial bifunctor of bidegree $(mp,mp)$ sending $(V,W)$ to  
$\Gamma^{m}(\gl(V^{(1)},W^{(1)}))$, where $\gl$ means $\Hom_k$.
The $\GL_n$-module $\Gamma^{m}(\gl_n^{(1)})$ is obtained by substituting $k^n$ for 
both $V$ and $W$ in $(V,W)\mapsto \Gamma^{m}(\gl^{(1)})(V,W)$.
Such substitution defines a functor $\Pol_d^d\to \Mod_{\GL_n}$ and for $n\geq d$ a
theorem of Franjou and Friedlander gives
$$\Ext^\bullet_{\Pol_d^d}(\Gamma^d\gl,F)\cong H^\bullet({\GL_n},F(k^n,k^n)).$$
The map from the left hand side to the right hand side goes by way of  $\Ext^\bullet_{\GL_n}(\Gamma^d\gl_n,F(k^n,k^n))$.
The invariant $\id^{\otimes d}\in \Gamma^d\gl_n$
gives a $\GL_n$-module map $k\to \Gamma^d\gl_n$ which allows one to go on to 
$\Ext^\bullet_{\GL_n}(k,F(k^n,k^n))=H^\bullet({\GL_n},F(k^n,k^n))$.

\subsection{The collapsing theorem for bifunctors}\label{bif collapse}
In order to explain the collapsing theorem for strict polynomial bifunctors in \cite{Touze 2013} we need to 
introduce a few counterparts of definitions  given above for strict polynomial functors.

So let $B\in \Pol^d_e$ be a bifunctor and let $Z\in \Gamma^e\V$. Then the parametrized bifunctor $B_Z$ is defined by 
parametrizing the covariant variable: $B_Z(V,W)=B(V,Z\otimes_k^{\Gamma^d} W)$.
Now assume $Z$ comes with an action of $\Gm$. Then $\Gm$ acts on $B_Z$ and we write $B^t_Z$ for the weight space of weight $t$.
The example we have in mind is $Z=E_r$ as in \ref{formality section}.

Define a Frobenius twist $B^{(r)}$ of $B$ by precomposition with $I^{(r)}$ in both variables:
$B^{(r)}(V,W)=B(V^{(r)},W^{(r)})$.

\begin{Theorem}[Collapsing of the Twisting Spectral Sequence \cite{Touze 2013}]
 Let $B\in \Pol^d_d$. There is a first quadrant spectral sequence 
$$E_2^{st}=\Ext_{\Pol^d_d}^s(\Gamma^d\gl,B^t_{E_r})\Rightarrow \Ext_{\Pol^{pd}_{pd}}^{s+t}(\Gamma^{dp}\gl,B^{(r)})$$
and this spectral sequence collapses.
\end{Theorem}

The proof of this theorem uses the themes that we seen above for ordinary polynomial functors:

\vbox{
\begin{itemize}
 \item adjoint of the twist,
 \item formality after twisting,
 \item untwisting a collapse.
\end{itemize}
}

Once one has the theorem one gets a much better grip on the connection between 
$ H^{2m}(\GL_{n},\Gamma^{m}(\gl_n^{(1)}))$ and
$ H^{2m}(\GL_{n},\bigotimes^{m}(\gl_n^{(1)}))$.
This then leads to the second generation construction of the classes of Touz\'e \cite{Touze 2013}.

\subsection{How the classes of Touz\'e help}

We find it hard to improve on the introduction to \cite{TvdK}.
Go read it.

\end{document}